\documentclass[11pt]{article}      
\usepackage{amsmath,amsthm,amsfonts,graphicx}
\usepackage{multirow}
\usepackage{multirow}
\usepackage{xcolor}
\def\smallddots{\mathinner{\raise7pt\hbox{.}\raise4pt\hbox{.}\raise1pt\hbox{.}}} 
\def\smallsdots{\mathinner{\raise1pt\hbox{.}\raise4pt\hbox{.}\raise7pt\hbox{.}}}

\DeclareMathOperator{\diag}{diag}

\DeclareMathOperator{\rank}{rank}

\DeclareMathOperator{\sign}{sign}

\numberwithin{equation}{section}
\numberwithin{table}{section}

\newtheorem{algorithm}{Algorithm}[section]
\newtheorem{example}{Example}[section]

\newtheorem{remark}{Remark}[section]

  \usepackage{enumitem}
\newlist{steps}{enumerate}{1}
\setlist[steps, 1]{label = Step \arabic*:}
\usepackage{booktabs}
\usepackage{subcaption}
\usepackage{tikz}
\usetikzlibrary{patterns, patterns.meta}
\usepgflibrary[patterns]

\begin{document}
 
\title{Superfast  1-Norm  Estimation\footnote{Dedicated to the memory of Nicholas J. Higham.}} 
\author{Victor Y. Pan} 
\author{Soo Go$^{[1],[a]}$ and Victor Y. Pan$^{[1, 2],[b]}$
\\ \\
$^{[1]}$ Ph.D. Programs in  Computer Science and Mathematics \\
The Graduate Center of the City University of New York \\
New York, NY 10036 USA \\
$^{[2]}$ Department of Computer Science \\
Lehman College of the City University of New York \\
Bronx, NY 10468 USA \\
$^{[a]}$
sgo@gradcenter.cuny.edu,~ 
$^{[b]}$
  victor.pan@lehman.cuny.edu \\ 
http://comet.lehman.cuny.edu/vpan/  \\
} 
\date{}

\maketitle


 
\begin{abstract}  
 \noindent
 A matrix algorithm is said to be  superfast (that is, runs at sublinear cost) if it involves much fewer  scalars and flops than an input matrix has entries. Such algorithms have been extensively studied and  widely applied in modern computations  for matrices with low displacement rank \cite{P01,P15, SW19} and more recently for low rank approximation of matrices, even though they fail on worst case inputs in the latter application.We extend this study to a new area. Our superfast algorithms consistently output  accurate 1-norm
estimates for real world   matrices, and we comment on some promising extensions of our surprisingly simple  techniques. With their further testing and refinement our algorithms can potentially be  adopted  in practical computations.
 
\end{abstract} 
 


\bigskip
\medskip


\section{Introduction}

 \subsection{Superfast matrix computations:   the State of the Art  (briefly)}\label{sstart}
 
Superfast (sublinear cost) algorithms are
highly important
in modern matrix computations,
e.g.,
where 
Big Data
are represented with matrices of immense size,  not fitting the primary memory of a computer. Quite
typically \cite{UT19},  such matrices can be closely 
approximated
with                                                                                  matrices of a low rank, 
 with which one can operate 
 superfast. Can we, however, compute {\em Low Rank Approximation (LRA)} of a matrix superfast?

Unfortunately, any superfast LRA algorithm fails on worst case inputs
 and even on the small matrix families of  our Appendix \ref{shrdin}, but superfast randomized  algorithms
of \cite{MW17,BW18,KLMMS23,LX23,CETW23},
devised for Symmetric Positive
Semidefinite  matrices or  other specified  matrix classes, output LRAs  whose  expected  Frobenius error norm is proved to be close to
  optimal. 

Even more relevant to our present work is superfast computation of  
LRA by means of    {\em Cross-Appro\-xi\-ma\-tion  (C-A) iterations},
which   apply the ADI method \cite{S16}   to LRA (see \cite{ALS24,OZ18},  and the references therein).
 For over two decades of their   
 worldwide application
C-A iterations have been consistently outputting
 accurate LRAs for a large and important class of matrices, even though
 adequate formal support for such  empirical 
performance is still a challenge.

The  algorithms of \cite{GLPa} 
enhance the power of C-A iterations for LRA
by combining them with compression of rank-$q$ approximation of a matrix into its rank-$r$ approximation for $q>r$ (cf. \cite{TYUC17,TYUC19}) and the techniques
of iterative refinement.
The algorithms of \cite{GLPa}
perform similarly where
 C-A iterations
are  replaced with superfast variation of
  the {\em generalized Nystr{\"o}m} LRA algorithm \cite{N20}.
In the latter algorithm 
an input matrix $M$ is pre- and post-multiplied by skinny  {\em Gaussian random matrices},
filled with
independent standard  Gaussian (normal) random variables, but \cite{GLPa}  make that stage and hence the entire algorithm superfast by replacing Gaussian multipliers with certain Ultrasparse
  matrices,
  said to be
 {\em 3-Abridged SRHT} matrices
 and specified in Appendix \ref{spreprmlt}.
 
 In  extensive numerical tests in  \cite{GLPa}, these superfast  algorithms of
both kinds 
  have  consistently computed 
  accurate  LRAs, although adequate formal support for that performance remains a challenge, as for C-A iterations. 

\subsection{Our results}
\label{soutl}

We extend the latter progress to  the estimation of the
1-norm of a matrix $M$ and hence
 its $\infty$-norm, which is the 1-norm of
the transpose $M^T$.

 LAPACK's 1-norm estimator essentially amounts to  at most 11 
(and on the average  less than 5)  multiplications
of $M$ and $M^T$ by vectors.
Such multiplications  and hence  LAPACK's  estimator
are superfast
where they are applied to a low rank matrix,
and so superfast LRA algorithms of \cite{GLPa} are immediately extended to superfast norm-1 estimators.

We tested such an extension
 by applying   LAPACK's estimator to the outputs of a most primitive algorithm of 
\cite{GLPa}, which uses just
3-Abridged SRHT multipliers and applies no iterative refinement or compression techniques. 

In our  numerical tests
 the relative output errors  ranged from 0.5 to 2.5, depending on an input matrix class.  Although quite reasonable, such output accuracy
is significantly inferior to LAPACK's.
We hoped
to improve
the accuracy
of our estimator
 by strengthening its LRA stage, but this has not worked in our tests where
we shifted back from \cite{GLPa}
to  \cite{N20} by replacing 3-Abridged SRHT  multipliers
with Gaussian
(see  Table
\ref{tabalgnrmestp1}). 
Since LRA of \cite{N20} is
expected to be near-optimal
under the Frobenius norm,
we had to
seek another way to superfast 1-norm estimator that would match the accuracy level of  LAPACK.

We succeeded when we incorporated {\em randomized  sparsification} into LAPACK's 1-norm estimator. 
Namely, before 
multiplying  $M$ or $M^T$ by a  vector,
we  sparsified it by keeping intact a
small fraction of its coordinates, chosen at random, and
 replacing all other coordinates with 0s.
 
We tested this estimator 
 for
synthetic and real world matrices $M$
in two ways --
 with and without  scaling the
sparse vectors, 
and in both cases our output  estimates for 1-norm $||M||_1$  have been consistently as accurate as LAPACK's.
This suggests that our estimators only fail on a narrow class of coherent 
inputs that are rare guests 
in Modern Matrix Computations.

Formal and empirical study of the class of such rare guests is a natural  challenge, possibly as hard as similar challenge for C-A iterations,
but with  further refinement 
and testing
our algorithms can potentially  be adopted in practice.
We  feel that 
highly efficient
matrix algorithms
supported  empirically
should be more
 widely used,
as this is now routine, e.g.,  for the algorithms in
Modern Cryptography
\cite{KL25}.

\subsection{Some
 extensions}
\label{sexts} 
 
 We succeeded based on  surprisingly simple algorithmic novelties and empirical feedback, and this should  motivate new  effort for devising
 superfast algorithms 
 for  matrix computations, possibly by means of extending the techniques of this paper and \cite{GLPa}. 
  
In Sec. \ref{scnc} 
(Conclusions) we briefly discuss
extension of LRA of 
\cite{GLPa} 
to {\em superfast least squares approximation}
(by following \cite{ALS24})
 and
of our randomized  sparsification techniques to
{\em superfast $p$-norm estimation}
 for any $p$. We  extend LRA of \cite{GLPa} to the latter problem in 
 Sec. \ref{sbckmnrm}.

\subsection{Related work}
Our work
exploits and extends impressive
power of LAPACK's 1-norm estimator, covered in
 \cite{ABB99,S98,
H02,HT00}, and the bibliography therein together with other customary matrix norm estimators. 
We were  largely motivated by 
the admirable recent progress in superfast LRA, cited 
in Sec. \ref{sstart}, 
 but to the best of our knowledge, no work on superfast matrix norm estimators as well as no work using our techniques for devising such estimators have ever been published yet.
Abridged SRHT matrices have been used in
\cite{PLSZ16,
PLSZ17,PLSZa}.
 
\subsection{Organization of the paper} 
In the rest of this section we briefly recall LAPACK’s 1-norm estimator.
 In Sec. \ref{srcp}  we extend it to superfast 1-norm estimator, which we  strengthen a little in  Sec. \ref{srci}  by means of incorporation
 of C-A iterations.  We  estimate
$p$-norms 
based on LRA in Sec. \ref{s2nrm},
 cover our numerical experiments in
Sec. \ref{ststs},
and devote  Sec. \ref{scnc} to conclusions.
 Every superfast norm estimator fails on small families of matrices of Appendix \ref{shrdin}.
In Appendix \ref{sspfcrd} we briefly recall random 
sketching for LRA and define Abridged SRHT matrices.
In Appendix \ref{sscl} we 
modify our scaling
 of Sec. \ref{srcp}.
    
\subsection{Background for 1-norm estimation}\label{sbcgour}
 
For a vector
${\bf v}=(v_i)^{n}_{i=1}\in\mathbb R^n$ and a
matrix $M=(m_{i,j})_{i,j=1}^{m,n}\in\mathbb R^{m\times n}$ 
   write
\begin{equation}\label{eqnrms}
||{\bf v}||_{\infty}:=
\max_{i=1}^n|v_i|,~||{\bf v}||_1:=\sum_{i=1}^n|v_i|,~
||M^T||_{\infty}:=||M||_1:=\max_{j=1}^n\sum_{i=1}^m|m_{i,j}|.
\end{equation}

 For a fixed column index $j$, $1\le j\le n$, and
 ${\bf e}_j$
 denoting the $j$th coordinate vector, LAPACK's 1-norm estimator (specified in \cite{ABB99}, \cite[Sec. 5.3.1]{S98},
 \cite[Sec. 15.2]{H02}) 
 recursively computes  the following vectors and
 vector norms:
\begin{equation}\label{equwx} 
{\bf u}:=(u_i)_{i=1}^{m}:=M{\bf e}_j,~{\bf w}:=(\sign(u_i))_{i=1}^{m},~{\bf x}:=M^T{\bf w},~||{\bf u}||_1={\bf w}^T{\bf u},~{\rm and}~||{\bf x}||_{\infty}.
\end{equation} 
The complexity  of this computation is  dominated at the stage of pre-multiplication of $M$ by a vector ${\bf w}^T$  filled with $\pm 1$.
For such a vector ${\bf w}$
and any column index $j$, $1\le j\le n$, obtain
\begin{equation}\label{eqMx} 
 ||M{\bf e}_j||_1\ge|{\bf w}^TM{\bf e}_j|=|{\bf x}^T{\bf e}_j|.
\end{equation}
Hence $||M{\bf e}_j||_1>||{\bf u}||_1$ if
\begin{equation}\label{eqMux} 
~ |{\bf x}^T{\bf e}_j|=||{{\bf x}}||_{\infty}>||{\bf u}||_1,
\end{equation}
and if ${\bf u}$ 
is a current candidate
 for maximizing the 1-norm of the column vectors of $M$, then  LAPACK updates ${\bf u}$ 
 with $M{\bf e}_j$, strictly increasing  the 1-norm.
Therefore,  in at most $n$ iterations for an $m\times n$  matrix $M$ a local maximum $||M{\bf e}_j||_1$ is reached. 

Under
this policy the  candidate 
value $||{\bf u}||_1$ is not  updated
unless 
(\ref{eqMux})
  holds, even where
 $||M{\bf e}_j||_1= ||M||_1>||{\bf u}||_1=|{\bf x}^T{\bf e}_j|$. Because of such   missing opportunities,
   one can end far  from the global maximum $||M||_1$ without noticing this, 
but  LAPACK   counters such mishaps by comparing $||M{\bf e}_j||_1$ for the initial and final vectors ${\bf e}_j$
with $||M{\bf g}||_1$
and  $||M{\bf h}||_1$, respectively,\footnote{The vector ${\bf h}$ was introduced  by N.J. Higham in 1988.} where
 \begin{equation}\label{eqlpck1nrm}
{\bf g}:=\left(\frac{1}{n},\dots,\frac{1}{n}\right)^T~{\rm and}~ {\bf h}:=\left(1,-1-\frac{1}{n-1},1+\frac{2}{n-1},
-1-\frac{3}{n-1},\dots,(-1)^{n-1}2\right)^T.
\end{equation}

\section{1-norm estimators with randomized sparsification}\label{srcp}

\subsection{Outline}\label{soutl}
To simplify our exposition we assume that $M$ is an $n\times n$  matrix and comment on extensions to  rectangular matrices $M$ in Remark \ref{remn}.
Hereafter we   write $a\ll b$ and $b\gg a$, for real $a$ and $b$,  to show that $a$ is much smaller than $b$ in context, which is more realistic than the asymptotic notation  $a=o(b)$,  hiding  overhead constants.   

Due to sparsification, our  superfast heuristic variations of LAPACK's estimator    
 involve much fewer than $n^2$ scalars  and  flops overall, although not fewer matrix-by-vector multiplications, but become  more likely
to end far from the norm $||M||_1$. In Secs.
\ref{ssprw/o} and \ref{sprem} we 
  modify our algorithms to counter this problem.

$k$-{\bf Sparsification}:                                                                                                                                                                                                                                                                                                                                                                                                                                                                                                                                                                                                                                                                                                                                                                                                                                                                                                                                                                                                                                                                                                                                                                                                                                                                                                                                                                                                                                                                                                                                                                                                                                                                                                                                                                                                                                                                                                                                                                                                                                                                                                                                                                                                                                                                                                                                                                                                                                                                                                                                                                                                                                                                                                                                                                                                                                                                                                                                                                                                                                                                                                                                                                                                                                                                                                                                                                                                                                                                                                                                                                                                                                                                                                                                                                                                                                                                                                                                                                                                                                                                                                                                                                                                                                       For two  integers $k$ and $n$, such that $1\le k\ll n$, \footnote{In our tests with $1024 \times 1024$ matrices we let $k$ take on the values 1, $3=\lfloor\log_2(\log_2 (n))\rfloor$, and $10=\log_2(n)$.} 
and for a vector                                                                                                                                                                                                                                                                                                                                                                                                                                                                                                                                                                                                                                                                                                                                                                                                                                                                                                                                                                                                                                                                                                                                                                                                                                                                                                                                                                                                                                                                                                                                                                                                                                                                                                                                                                                                                                                                                                                                                                                                                                                                                                                                                                                                                                                                                                                                                                                                                                                                                                                                                                                                                                                                                                                                                                                                                                                                                                                                                                                                                                                                                                                                                                                                                                                                                                                                                                                                                                                                                                                                                                                                                                                                                                                                                                                                                                                                                                                                                                                                                                                                                                                                                                                                       ${\bf v}$ of dimension $n$, do
\begin{enumerate}
\item 
 keep the coordinates of the vector  at $k$ random positions chosen under the uniform probability distribution  on the set $\{1,2,\dots,n\},$   
 \item  then set the $n-k$
remaining coordinates  to 0 and output the resulting vector   $\bar {\bf v}$. 
\end{enumerate}

We first compute the
$k$-sparsified  vectors 
$\bar{\bf g}$
and $\bar{\bf h}$
and then  multiply $M$ by the normalized
output vectors
$\widehat{\bf g}:=\bar{\bf g}/||\bar{\bf g}||_1$ and 
$\widehat{\bf h}:=\bar{\bf h}/||\bar{\bf h}||_1$. Unlike LAPACK, we do this 
 at the initialization stage and always  
 approximate 
$||M||_1$ with  $||M{\bf e}_j||_1$ for some index $j$ -- never with $||M\widehat{\bf g}|| _1$ or $||M\widehat{\bf h}||_1$. 

We modify LAPACK's policy
of updating the candidate  index $j$ in two ways. 
In the next subsection
we  update it wherever
this increases
the 1-norm $||M{\bf e}_j||_1$; we test this provision superfast.
In Sec. \ref{sprem} we
 also require for updating that $||M{\bf e}_j||_1\ge \alpha ||{\bf x}||_{\infty}$ for ${\bf x}$ of (\ref{equwx}) and a fixed $\alpha\ge 1$.

In Sec. \ref{sn/k} we cover our tests  for $\alpha$ equal to the {\em sparsification ratio} $n/k$ -- to compensate 
for the decrease of
$||{\bf x}||_1$
due to $k$-sparsification
of vector ${\bf w}$ of dimension $n$. In Appendix \ref{sscl} we
 recursively adapt $\alpha$ according to a feedback from observed impact
on the  norm $||M{\bf e}_j||_1$. 

Next we specify our modifications of LAPACK's estimator. 

\subsection{A sparsified 1-Norm Estimator without scaling}\label{ssprw/o}

\begin{algorithm}\label{algnormest2}
{\rm [A sparsified 1-Norm Estimator 1.]}
 

\begin{description}


\item[{\sc Input:}] two positive integers $k$ and $n$ such that $k\ll n$, 
a real $n\times n$ matrix $M$, two vectors
${\bf g}$ and ${\bf h}$ of (\ref{eqlpck1nrm}),
and a  tolerance TOL $\ge 2$.


\item[{\sc Output:}] A column index  $j$, $1\le j\le n$, and a value $\nu\ge 0$ [expected to  approximate  $||M||_1$].


\item[{\sc Initialization:}]
\begin{enumerate}
\item 
Compute the vectors $\bar {\bf g}$ and $\bar {\bf h}$, denoting  the  $k$-sparsified
vectors ${\bf g}$ and $ {\bf h}$. \item  Compute sparse normalized vectors
$\widehat{\bf g}:=\bar{\bf g}/||\bar{\bf g}||_1$ and $\widehat {\bf h}:=\bar{\bf h}/||\bar{\bf h}||_1$.
\item
Compute
the vectors 
$M\widehat {\bf g}$ and $M\widehat{\bf h}$ and their 1-norms. 

\item
Write ${\bf u}  = (u_i)^n_{i=1} := M\widehat {\bf g}$ if $||M\widehat {\bf g}||_1 \ge ||M\widehat {\bf h}||_1$; 
write  ${\bf u}= (u_i)^n_{i=1} := M\widehat {\bf h}$ otherwise.

\item 
Write $\gamma:=0$
and $\nu_0:=-1$.
\end{enumerate}


\item[{\sc Computations:}]
\begin{enumerate}

\item  
Compute the vectors $ {\bf w}:= (\sign(u_i))_{i=1}^n$, $\bar{\bf w}$
(the $k$-sparsified vector  $ {\bf w}$), and $ {\bf x}:= (x_j)_{j=1}^n= M^T  \bar{\bf w}$ and increase the value $\gamma$ by 1.

\item Choose an index $j(\gamma)$ such that $|x_{j(\gamma)}| = || {\bf x}||_\infty$.
\item
Write  ${\bf u} := M {\bf e}_{j(\gamma)}$, compute $\nu_{\gamma}:=||M {\bf e}_{j(\gamma)}||_1 $.
\item
Stop and  output $j:=j(\gamma-1)$ 
and $\nu:=\nu_{\gamma-1}$ if  $\nu_{\gamma-1}\ge  \nu_{\gamma}$.
\item
Stop and  output $j:=j(\gamma)$ and $\nu:=
\nu_{\gamma}$ 
if $\gamma = \text{TOL}$. 
\item
Go to Stage 1.

\end{enumerate}

\end{description}
\end{algorithm}

\begin{remark}\label{reghj}
(i) The initial vector ${\bf u}$  is only used at Stage 1 of the computations 
(for $\gamma=0$)  and then is changed at Stage 3. Hence our estimators, unlike 
LAPACK's, never output the 1-norms
 $||M{\bf g}||_1/||{\bf g}||_1$
 or $||M{\bf h}||_1/||{\bf h}||_1$. (ii) Computations of Alg. \ref{algnormest} 
cannot stop at Stage 4 unless $\gamma>1$
because $\gamma_0=-1<
\min\{||{\bf x}||_\infty, \nu_{1}\}$.
\end{remark}

\subsection{Computational cost estimates}
Our
computational cost estimates
{\em omit the cost of the generation of random parameters,}
assumed to be dominated. Every multiplication of $M$ or $M^T$ by a $k$-sparsified vector involves at most $kn$ entries of $M$ and  $kn$ flops. Overall  the initialization
stage and $s$ stages of computations 
 involve $(s+2)kn$ entries of $M$, $(s+2)kn$ flops and in addition 
$(n+2)s$
pairwise comparisons of real numbers,  $(n-1)(s+3)$ additions, $2n$ divisions, and computation of   the signs of $ks$ real numbers. 
 The number of entries of $M$ and all operations involved
is much smaller than $n^2$ for $s\ll n/k$.

\subsection{Sparsified 1-norm estimation with scaling}\label{sprem}
 
Our second 1-norm estimator
includes a fixed $\alpha\ge 1$ into its inputs and otherwise 
differs from Alg. \ref{algnormest2}
only at Stage 4 of    computations, 
by involving  an extra pair of  comparison and multiplication.
 
\begin{algorithm}\label{algnormest}
{\rm [A sparsified LAPACK 1-Norm Estimator with scaling.]}
 

\begin{description}


\item[{\sc Input:}] a real $\alpha\ge 1$; otherwise as in Alg. \ref{algnormest2}.


\item[{\sc Output and  Initialization:}]
as in Alg. \ref{algnormest2}.

\item[{\sc Computations:}]
as in Alg. \ref{algnormest2}
except for Stage 4, changed as follows:

Stop and output $j:=j(\gamma)$ and $\nu:=\nu_{\gamma}$ if   $\nu_{\gamma-1} \ge \min\{\alpha||{\bf x}||_\infty, \nu_{\gamma}\}$. 

\end{description}
\end{algorithm}

\section{1-norm estimation incorporating  C-A  iterations}\label{srci}

\subsection{ C-A  iterations for rank-1 approximation}
\label{smxvll1}  

Next we   enhance our 1-norm estimators by means of incorporation   
of the {\em maxvol} algorithm
  of \cite{GOSTZ10},   which implements  C-A iterations.  
  We only use {\bf maxvol 1}, that is,  maxvol  for rank-1  approximation.

 \begin{figure} 
[ht]  
\centering
\centering
\begin{tikzpicture}

\def\rows{20}
\def\cols{16}
\def\cellsize{0.25}

\draw (0, 14*\cellsize) -- (\cols*\cellsize, 14*\cellsize);
\draw (0, 15*\cellsize) -- (\cols*\cellsize, 15*\cellsize);

\draw (3*\cellsize, 0) -- (3*\cellsize, \rows*\cellsize);
\draw (4*\cellsize, 0) -- (4*\cellsize, \rows*\cellsize);

\draw (12*\cellsize, 0) -- (12*\cellsize, \rows*\cellsize);
\draw (11*\cellsize, 0) -- (11*\cellsize, \rows*\cellsize);

\draw (3*\cellsize, 4*\cellsize) rectangle (4*\cellsize, 5*\cellsize);

\draw[pattern={Lines[angle=-45,distance={1.25pt}]}]
    (3*\cellsize, 14*\cellsize) rectangle (4*\cellsize, 15*\cellsize);
    
\draw[pattern={Lines[angle=45,distance={1.25pt}]}]
    (11*\cellsize, 14*\cellsize) rectangle (12*\cellsize, 15*\cellsize);

\draw[pattern={Lines[angle=90,distance={1.25pt}]}]
    (3*\cellsize, 4*\cellsize) rectangle (4*\cellsize, 5*\cellsize);

\draw[thick] (0,0) rectangle (\cols*\cellsize, \rows*\cellsize);

\end{tikzpicture}
\caption{The first three  steps of maxvol 1 output three entries marked with
strips.}
\label{fig4}
\end{figure}
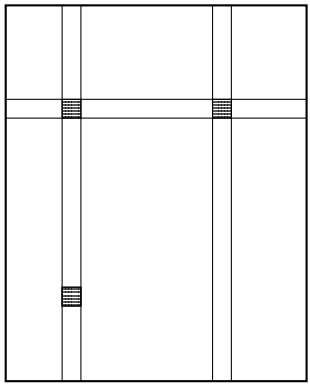  

\begin{algorithm}\label{algcrit}
{[Maxvol 1, see Fig. \ref{fig4}.]}


\begin{description}


\item[{\sc Input:}] 
a real $n\times n$ matrix $M=(m_{i,j})_{i,j=1}^{m,n}$.


\item[{\sc Output:}]
A pair of indexes $i$ and  $j$  such that 
\begin{equation}\label{eqij} 
|m_{i,  
 j}|\ge\max\{\max_{g=1}^m|m_{g,j}|, \max_{h=1}^n|m_{i,h}|\}.
\end{equation} 


\item[{\sc Initialization:}] \quad
 Fix a column index $j$,
 $1\le j\le n$, and 
  compute $i:= {\rm argmax}_{g=1}^m|m_{g,j}|$. 
 
 \item[{\sc Computations:}] \quad

Compute $\bar j:= {\rm argmax}_{h=1}^n|m_{ i,h}|$. 

Stop and output $i$
and $j$ if $\bar j=j$, that is, if 
$|m_{i,\bar j}|=|m_{i,j}|$.

Otherwise write  $j:=\bar j$ and compute $\bar i:= {\rm argmax}_{g=1}^m|m_{g,j}|$. 

Stop and output $i$
and $j$ if $\bar i=i$, that is, if 
$|m_{\bar i,j}|=|m_{i,j}|$.

Otherwise write  $i:=\bar i$ and repeat.

\end{description}
\end{algorithm} 

\begin{remark}\label{rerclmn}
 We can define {\rm dual maxvol 1} by applying  maxvol 1 to the transpose $M^T$.
\end{remark}

Every computation of
argmaximum requires $m-1$ or $n-1$ comparisons among the values $|m_{i,j}|$
 and  results  
in either  termination of computations or  a
strict increase
of  the current maximal $|m_{i,j}|$. Hence  the algorithm
stops at a local maximum 
 (\ref{eqij}) in at most $mn-1$  comparisons,
 but empirically it tends to run  superfast.
 This is impossible for worst case inputs
 and even for the small input
 families of Appendix \ref{shrdin},
 but empirically such hard inputs are rarely encountered.  In the case of a Gaussian random input matrix $M$,
maxvol 1
superfast  computes a close approximation of   
the global maximum 
\begin{equation}\label{eqij||}   
||M||_0:=|M|:=\max_{g,h=1}^{m,n}|m_{g,h}|
\end{equation} 
 with a high probability, according to formal and  empirical study in \cite{O17}.

\begin{remark}\label{reini}
{\rm[Alternative initialization of maxvol 1.]} Let us refer as {\bf Algorithm \ref{algcrit}-k} to the modification of maxvol 1 
where  initial column index is defined by the output of superfast Alg.    
\ref{algnormest2},
 unlike the standard initialization  recipe of maxvol 1.
In
the test results in  Table 
\ref{tabalg31-33comparison}
(for $n=1024$),  
Alg. \ref{algcrit}
superseded Alg.
\ref{algcrit}-k in accuracy only slightly and
 only for some input classes. One can similarly employ the output of 
Alg. \ref{algnormest}
instead of that of Alg.
\ref{algnormest2}.
\end{remark}

 \subsection{ Incorporations of maxvol 1 into 1-norm estimators}
\label{smxvll1}

 Notice that
maxvol 1 outputs a column  index $j$
such that
$m~||M{\bf e}_j||_1\ge ||M||_1$.
Hence these
$j$ and $||M{\bf e}_j||_1$
 can be a good   starting pair for our 
1-norm estimators, and we can
change
the entire Initialization
stage of both 
 Algs.   
\ref{algnormest2} and
\ref{algnormest} as follows:
\medskip

 For an index $j$ output by
maxvol 1,  fix ${\bf u}:=M{\bf e}_j                                                                                                                                                                                                                                                                                               $ and then 
write $\gamma:=0$
and $\nu_0:=-1$.                                                                                                                                                                                                                                                                                                                                                                                                                                                                                                                                                                                                                                                                                                                                                                                                                                                                                                                                                                                                                                                                                                                                                                                                                                                                                                                                                                                                                                                                                                                                                                                                                                                                                                                                                                                                                                                                                                                                                                                                                                                                                                                                                                                                                                                                                                                                                                                                                                                                                                                                                                                                                                                                                                                                                                                                                                                                                                                                                                                                                                                                                                                                                                                                                                                                                                                                                                                                                                                                                                                                                                                                                                                                                                                                                                                                                                                                                                                                                                                                                                                                                                                                                                                                                                              
\medskip
 
 In Figs. \ref{figalg32acc} and \ref{figalg32its} and 
Table~\ref{tabalg21-32comparison} we display our test 
results for a similar 1-norm estimator where we keep the Initialization  stage intact but replace 
 Stage 2 of computations of Alg.  
\ref{algnormest2}
 with new Stages 
2--4.~\footnote{One can similarly modify  Alg.  
\ref{algnormest}.} For convenience we next                                                                                                                                                  display 
 the entire modified algorithm. 
 In that display we allow tol invocations of  maxvol 1
  for a fixed tolerance tol$\ge 1$,
 but in our tests  we always  fixed tol=1 because
 already
in a single  invocation of maxvol 1
we obtain a  pair of indexes
 $i$  and $j$ such that $|m_{i,j}|\approx |M|$.

\begin{algorithm}\label{algnormest1}
{\rm [A sparsified LAPACK 1-Norm Estimator with  C-A steps.]}
 

\begin{description}


\item[{\sc Input:}] two positive integers $k$ and $n$ such that $k\ll n$, a real $n\times n$ matrix $M$, two vectors
${\bf g}$ and ${\bf h}$ of (\ref{eqlpck1nrm}),
and two  tolerance  values, {\rm tol} and  {\rm TOL}, $1\le {\rm tol}<{\rm TOL}$.


\item[{\sc Output:}]A column index  $j$, $1\le j\le n$, and a value $\nu\ge 0$ [expected to  approximate  $||M||_1$].


\item[{\sc Initialization:}]
\begin{enumerate}
\item 
Compute the vectors $\bar {\bf g}$ and $\bar {\bf h}$, denoting  the  $k$-sparsified
vectors ${\bf g}$ and $ {\bf h}$. \item  Compute sparse normalized vectors
$\widehat{\bf g}:=\bar{\bf g}/||\bar{\bf g}||_1$ and 
$\widehat{\bf h}:=\bar{\bf h}/||\bar{\bf h}||_1$.
\item
Compute
the vectors 
$M\widehat {\bf g}$ and $M\widehat{\bf h}$ and their 1-norms. 

\item
Write ${\bf u}  = (u_i)^n_{i=1} := M\widehat {\bf g}$ if $||M\widehat {\bf g}||_1 \ge ||M\widehat {\bf h}||_1$;
write ${\bf u}= (u_i)^n_{i=1} := M\widehat {\bf h}$ otherwise.

\item
Write $\gamma:=0$
and $\nu_0:=-1$. 
\end{enumerate}


\item[{\sc Computations:}]
\begin{enumerate}

\item  
Compute the vectors $ {\bf w}:= (\sign(u_i))_{i=1}^n$, $\bar{\bf w}$
(the $k$-sparsified vector  $ {\bf w}$), and $ {\bf x}:= (x_j)_{j=1}^n= M^T  \bar{\bf w}$   and increase the  value $\gamma$ by 1.
\item Choose an index $j(\gamma)$ such that $|x_{j(\gamma)}| = || {\bf x}||_\infty$.
\item If $\gamma \le {\rm tol}$, then apply Alg. \ref{algcrit} initialized with the index $j:=j(\gamma)$ to obtain a pair of row and column indexes $i^*$ and $j^*$. If $||M{\bf e}_{j*}||_1 > ||M{\bf e}_{j}||_1$, update $j(\gamma):=j^*$.
  \item
Write  ${\bf u} := M {\bf e}_{j(\gamma)}$ and compute $\nu_{\gamma}:=||M {\bf e}_{j(\gamma)}||_1 $.
\item
Stop and  output $j:=j(\gamma-1)$ 
and $\nu:=\nu_{\gamma-1}$ if  $\nu_{\gamma-1}\ge \nu_{\gamma}$.
\item
Stop and  output $j:=j(\gamma)$ and $\nu:=\nu_{\gamma}$ if $\gamma = \text{TOL}$.
 \item
Go to Stage 1.
\end{enumerate}
\end{description}
\end{algorithm}

 \begin{remark}\label{remn}
To extend  Algs.
\ref{algnormest2}, \ref{algnormest},
and \ref{algnormest1} to a rectangular 
 $m\times n$ matrix
 $M$, apply $k_n$-sparsification  procedure to
the vectors ${\bf g}$
 and ${\bf h}$
 of dimension $n$ at the initialization stage and apply $k_m$-sparsification  procedure to
  the vectors ${\bf w}$
 of dimension $m$ at Stage 1 of computations.
\end{remark} 

                                                                                                                                                                                                                                                                                                                                                                                                                                                          
\section{Application of LRA to matrix norms estimation} \label{s2nrm}

\subsection{Background on LRA}\label{sbcgrlra}

\begin{figure}[ht]
\centering
\begin{tikzpicture}[scale=1]

\def\gap{0.3}

\draw (-3.5,0) rectangle (-1, 3);
\node at (-2.1,1.5) {$M$};

\node at (-0.5,1.5) {$=$};

\draw (0,0) rectangle (0.5,3);
\node at (0.25,1.5) {$X$};

\draw (0.5+\gap,2.5) rectangle (1+\gap,3);
\node at (0.75+\gap,2.75) {$Y$};

\draw (1+2*\gap,2.5) rectangle (3.5+2*\gap,3);
\node at (2.25+2*\gap,2.75) {$Z$};

\node at (4+2*\gap,1.5) {$+$};

\draw (4.5+2*\gap,0) rectangle (7+2*\gap, 3);
\node at (5.75+2*\gap,1.5) {$E$};

\end{tikzpicture}
\caption{LRA  of $M$.}
\label{figlra}
\end{figure}
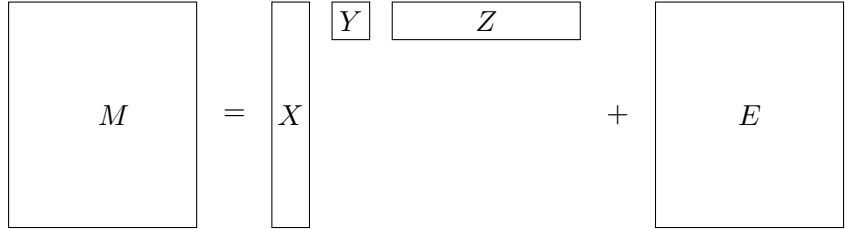

A triplet  
$\{X,Y,Z\}$
of matrices is said to be 
a {\em rank-$r$ approximation}\footnote{Here and hereafter ``rank-$r$" means ``rank of at most $r$".}  of an $m\times n$  matrix $M$ (see Fig. \ref{figlra}) 
 if
\begin{equation}\label{eqXYZkl}
X\in \mathbb C^{m\times k},~ Y\in \mathbb C^{k\times \ell},~ Z\in \mathbb C^{\ell\times n},~r\le\min\{k,\ell\},~{\rm and}~E:=
M-~XYZ\approx 0~{\rm in~ context}.
\end{equation}
 
  For $k=\ell=r$ and the    
identity matrix  $Y$,
the triplet degenerates into a pair $\{X,Z\}$,
but we can also readily compress  
a triplet $\{X,Y,Z\}$ into the pair $\{XY,Z\}$
of 
 $m\times \ell$ and $\ell\times n$ matrices or
 the pair $\{X,YZ\}$ 
of  $m\times k$ and $k\times n$
matrices.

 $\sigma_j(M)\ge 0$
 denotes the $j$th largest singular value of a matrix $M$.
The matrix $M^{(r)}$,
for $r\le {\min\{m,n\}}$, denotes the $r$-{\em truncation} 
obtained by means of setting to 0 all  singular values of  $M$ but its $r$ top (largest) ones. The compact SVD of $M^{(r)}$,   aka the $r$-{\em top SVD} of $M$, is  
  a  rank-$r$ approximation  of $M$
  given by the triplet $\{X,Y,Z\}=\{U^{(r)},\Sigma^{(r)},V^{(r)}\}$ such that $M^{(r)}=U^{(r)}\Sigma^{(r)}V^{(r)T}$, $\Sigma^{(r)}:=\diag(\sigma_j(M))_{j=1}^r$, and $U^{(r)}\in \mathbb R^{m\times r}$ and $V^{(r)}\in \mathbb R^{n\times r}$  are  two matrices  filled with
 the associated  orthonormal
 left and right singular vectors of $M$. This rank-$r$ approximation  of $M$ is optimal under the spectral and Frobenius matrix norms $||\cdot||_2$ and $||\cdot||_F$, respectively, by virtue of Eckart-Young-Mirsky theorem.   
 
 \subsection{Background for matrix norm estimation}\label{sbcknrm}
 
In addition to (\ref{eqnrms})
and (\ref{eqij||}),
define the $p$-norm  
\begin{equation}\label{eqpnrm}
||M||_p:=\max_{||{\bf x}||_p=1}||M{\bf x}||_p,~~||M||_2=\sigma_1(M)
\end{equation} 
for a matrix
$M=(m_{i,j})_{i,j=1}^{m,n}\in \mathbb C^{m\times n}$, vectors ${\bf x}\in \mathbb C^{n\times 1}$,
and $p\ge 1$. Recall  that 
\begin{equation}\label{eqnrm2p}
\frac{1}{n^{|1/p-1/2|}}\le\frac{||M||_p}{||M||_2}\le 
n^{|1/p-1/2|}~{\rm for}~m=n,  
\end{equation}
\begin{equation}\label{eqnrm02}
||M||_0\le||M||_2\le \sqrt{mn}~||M||_0~{\rm for}~||M||_0:=|M|~{\rm of}~(\ref{eqij||})
\end{equation}
(see  \cite[Sec. 7.3]{H02} and \cite[Sec. 2.3.2]{GL13}, respectively)
and  that the  known $p$-norm estimators are reduced
 to 
recursive multiplication of $M$ and $M^T$
by vectors
 (see  \cite{B74,D83}, \cite[Secs. 15.2 and 15.4]{H02},
and \cite[Sec. 5.3.4]{S98}).

 \subsection{Matrix norm estimation via LRA}\label{sbckmnrm}

Alg. \ref{algnrmest} represents 
our LRA-based $p$-norm estimator.

\begin{figure}[h!]
\begin{algorithm}\label{algnrmest}
{\rm [A  Matrix Norm Estimator via LRA.]}
 

\begin{description}


\item[{\sc Input:}] 
 $M\in \mathbb C^{m\times n}$
and $p\ge 1$  or $p=0$.

\item[{\sc Output:}] $\nu_p\approx ||M||_p$.

\item[{\sc Initialization:}]
Fix an integer $r>1$.

\item[{\sc Computations:}]
\begin{enumerate}
\item \quad 
Compute a close
 rank-$r$
approximation $\{X,Y,Z\}$ of $M$ under the 2-norm. 
\item
 Compute $\nu_p\approx ||XYZ||_p$ by  applying a known $p$-norm estimator.

\end{enumerate}
\end{description}
\end{algorithm} 

\end{figure}

Alg. \ref{algnrmest} is superfast 
 if (i) so is the LRA 
algorithm  
applied at its Stage 1, e.g., if this is an algorithm of  \cite{GLPa}, (ii)
a known $p$-norm estimator,
applied  at Stage 2,
converges in  a reasonably bounded number of multiplications of $M$ and $M^T$ by vectors, and (iii)  $\rank(XYZ)$ is small, so that the matrix $XYZ$ is multiplied by a vector superfast.
 
To estimate the output error
$|\nu_p-||M||_p|$ of a $p$-norm estimate $\nu_p\approx ||M||_p$, 
 apply triangle inequality and obtain
$$|\nu_p-||M||_p|\le |\nu_p-||XYZ||_p|+\Delta_p~{\rm for}~\Delta_p:=||M-XYZ||_p.$$

Assume
that the  error bounds 
$\Delta_2:=||M-XYZ||_2$  and
$|\nu_p-||XYZ||_p|$ are small 
and extend the  norm bound $\Delta_2$ to  $\Delta_p$ for all $p$ by combining 
 Eqns. \ref{eqij||}, \ref{eqnrm2p}
and \ref{eqnrm02}.

\begin{remark}\label{refrb}
 Extension of Alg. \ref{algnrmest}
to estimation of the Frobenius norm of  $M$ is immediate. To extend
the output error
bound, apply 
\cite[Eqn. 2.3.7]{GL13}.

\end{remark}

\section{Numerical experiments for 1-norm estimation}\label{ststs}

 We implemented our algorithms   in Python,  run all tests  on 64-bit MacOS 14.5 (ARM64) with 16GB Memory, and used numpy and scipy packages for the numerical linear algebra routines.

\subsection{Input Matrices} \label{sec_input}

\noindent {\bf Real-world Input Matrices.} 
 
The matrix {\bf Shaw}  
is from a one-dimensional image restoration model problem. 
 
The matrix  {\bf Gravity}   
is
from a one-dimensional gravity surveying model problem. 

Both Shaw and  Gravity 
 are $1000\times 1000$ dense real matrices  
having low numerical rank; we pad them with 0s  to increase their size to $1024\times 1024$. 
They represent discretization of Integral Equations, provided in the built-in problems of the Regularization Tools\footnote{See 
 http://www.math.sjsu.edu/singular/matrices and 
  http://www.imm.dtu.dk/$\sim$pcha/Regutools
  
 For more details see Ch. 4 of 
  http://www.imm.dtu.dk/$\sim$pcha/Regutools/RTv4manual.pdf }. 

 Our third 
input matrix is from the discretization of a single layer potential ({\bf SLP}) operator (see   \cite[Sec. 7.1]{HMT11} for  details).
It has size  $1024\times 1024$.  
\medskip

{\bf Synthetic Input Matrices:} We generated 
synthetic $1024\times 1024$ input  matrices of 
six classes,
the first four of them  are the products $U\Sigma V^T$, 
for $\Sigma = \textrm{diag}(\sigma_i)_{i=1}^{1024}$ and $U$ and $V$ being the orthogonal matrices of the left and right singular vectors of a 
$1024\times 1024$  Gaussian random matrix, respectively.

To generate {\bf Fast Decay} matrices, i.e., those having fast decaying spectra, we let
  $\sigma_i = 1$ for $i = 1, 2, \dots, 20$, 
$\sigma_i = \frac{1}{2^{i-20}}$ for $i = 21, \dots, 100$, and $\sigma_i = 0$ for $i > 100$.

To generate {\bf Slow Decay} matrices, i.e., those having slowly decaying spectra, we let 
 $\sigma_i = 1$ for $i = 1, 2, \dots,20$, 
and $\sigma_i = \frac{1}{(1 + i - 20)^2}$ for $i > 20$. 

The synthetic input matrices of the next three types are from \cite{H88}.

 To generate  matrices with a {\bf Single Small Singular Value}, we let $\sigma_i = 1$ for $i = 1, \dots, 1023$ and $\sigma_{1024} \in [10^{-16}, 10^{-3}]$.

 To  generate matrices with a {\bf Single Large Singular Value}, we let  $\sigma_1 \in [10^3, 10^{16}]$ and 
 $\sigma_i = 1$ for $i = 2, \dots, 1024$.
 
To generate  {\bf Cauchy} matrices, we first  randomly defined two vectors ${\bf e}$ and ${\bf f}$ of dimension $1024$ with independent random coordinates chosen  uniformly  from the interval  $[0, 1]$. Then we computed two vectors  
${\bf x} = a ${\bf 1}$+ (b-a){\bf e}$ and ${\bf y} =  c ${\bf 1}$+ (d-c){\bf f}$ for some fixed scalars $a, b, c$, and $d$ and vector ${\bf 1} =  (1)_{i=1}^{1024}$. Finally, we defined our input matrix  $M: = ((x_i - y_j)^{-1})_{i, j=1}^{1024}$. In our implementation, we fixed $a=0$, $b=100$, $c=100$, and $d=200$.

We generated  dense $1024\times 1024$  {\bf Random} matrices by filling their entries   with $-1$, 0, and 1 sampled with equal probability.

\subsection{Test Results for Algorithm~\ref{algnormest2}}\label{sn/k}
\begin{figure}[h!]
\centering
\includegraphics[width=.9\textwidth]{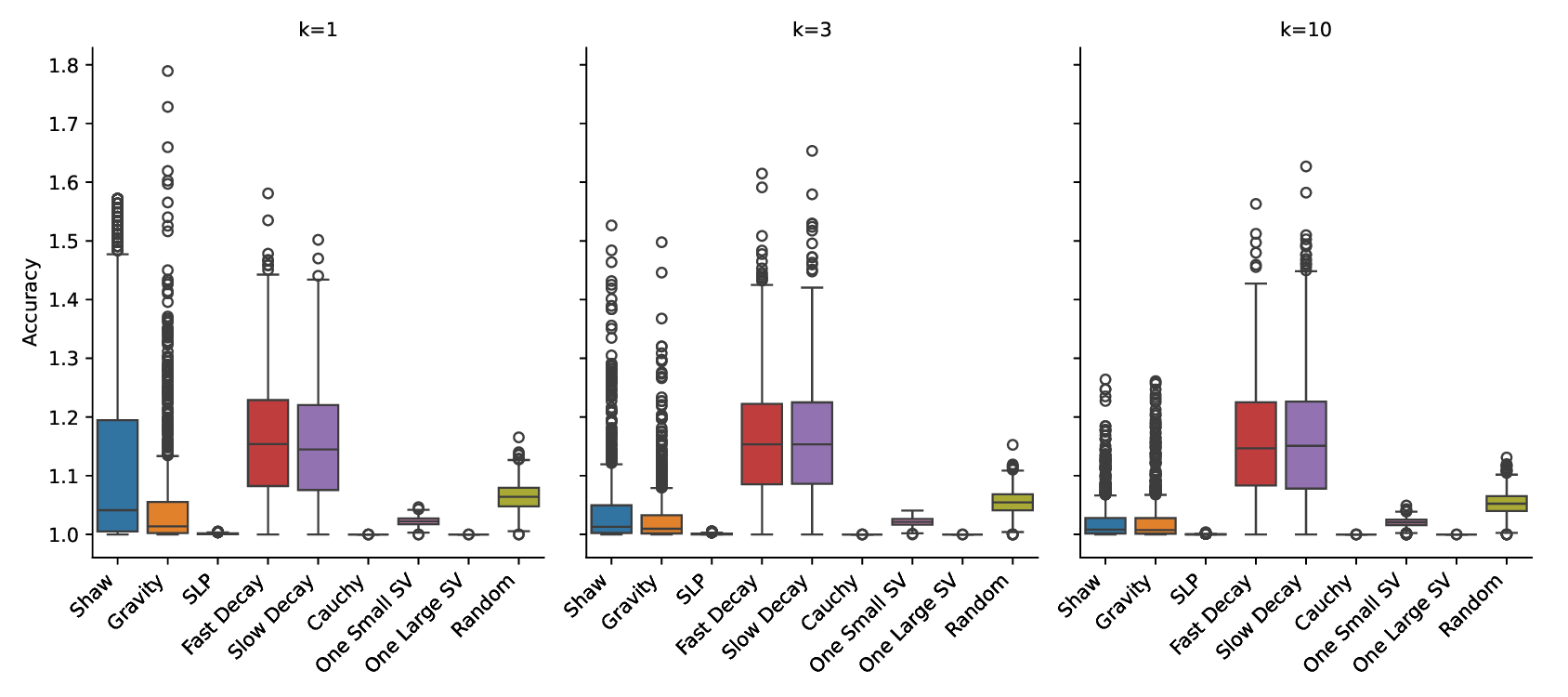}
\caption{Output errors of Alg.~\ref{algnormest2}.}
\label{figalg21acc}
\end{figure}

\begin{figure}[h!]
\centering
\includegraphics[width=.9\textwidth]{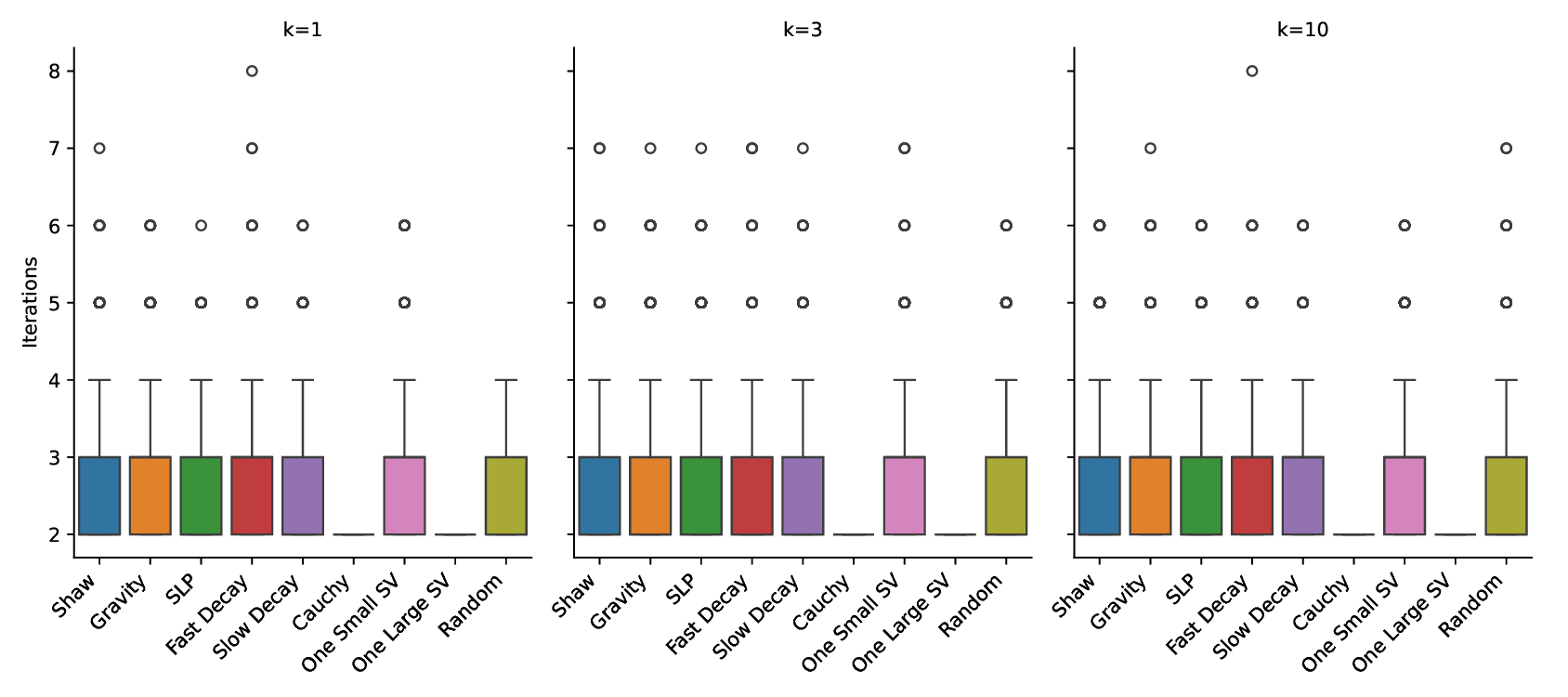}
\caption{Number of iterations of Alg.~\ref{algnormest2}.}
\label{figalg21its}
\end{figure}

We ran 1000 tests for Alg.~\ref{algnormest2} per every pair of a matrix type and an integer $k$,  $k = 1$, $k = \lfloor \log_2(\log_2(1024)) \rfloor=3$, and $k = \log_2(1024) = 10$. We fixed  TOL= 10, and measured the accuracy of the estimates by the ratio $\frac{||M||_1}{\text{est}}$ where est denoted the estimate of $||M||_1$ returned by the algorithm. The ratio took on its minimal value 1 where the estimate est was optimal. Figs. \ref{figalg21acc} and \ref{figalg21its} show  estimates for 1-norm of $M$  within a factor of 2 in  2, 3, or 4 iterations.

\subsection{Test Results for Algorithm~\ref{algnormest}}
\begin{figure}[h!]
\centering
\includegraphics[width=.9\textwidth]{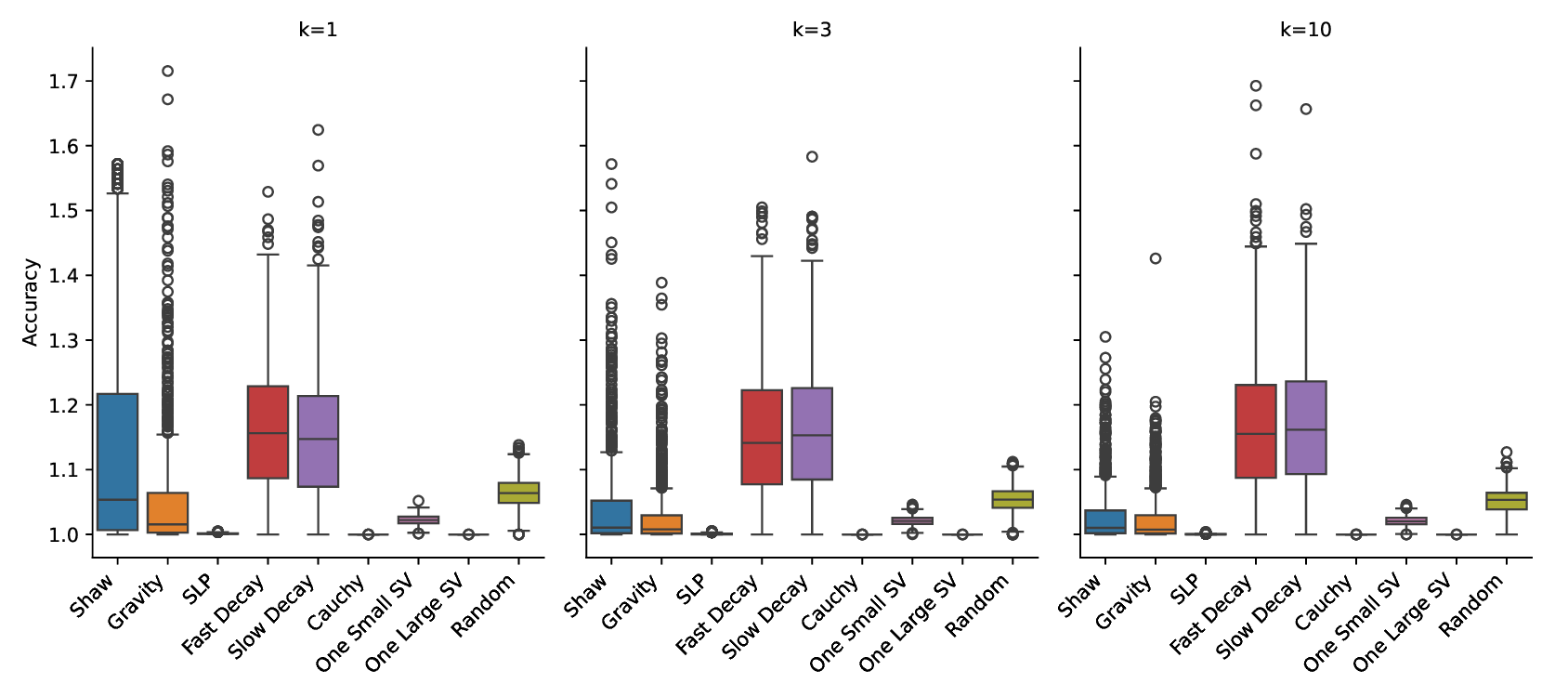}
\caption{Output errors of Alg.~\ref{algnormest}\label{figalg22acc}.}
\end{figure}

\begin{figure}[h!]
\centering
\includegraphics[width=.9\textwidth]{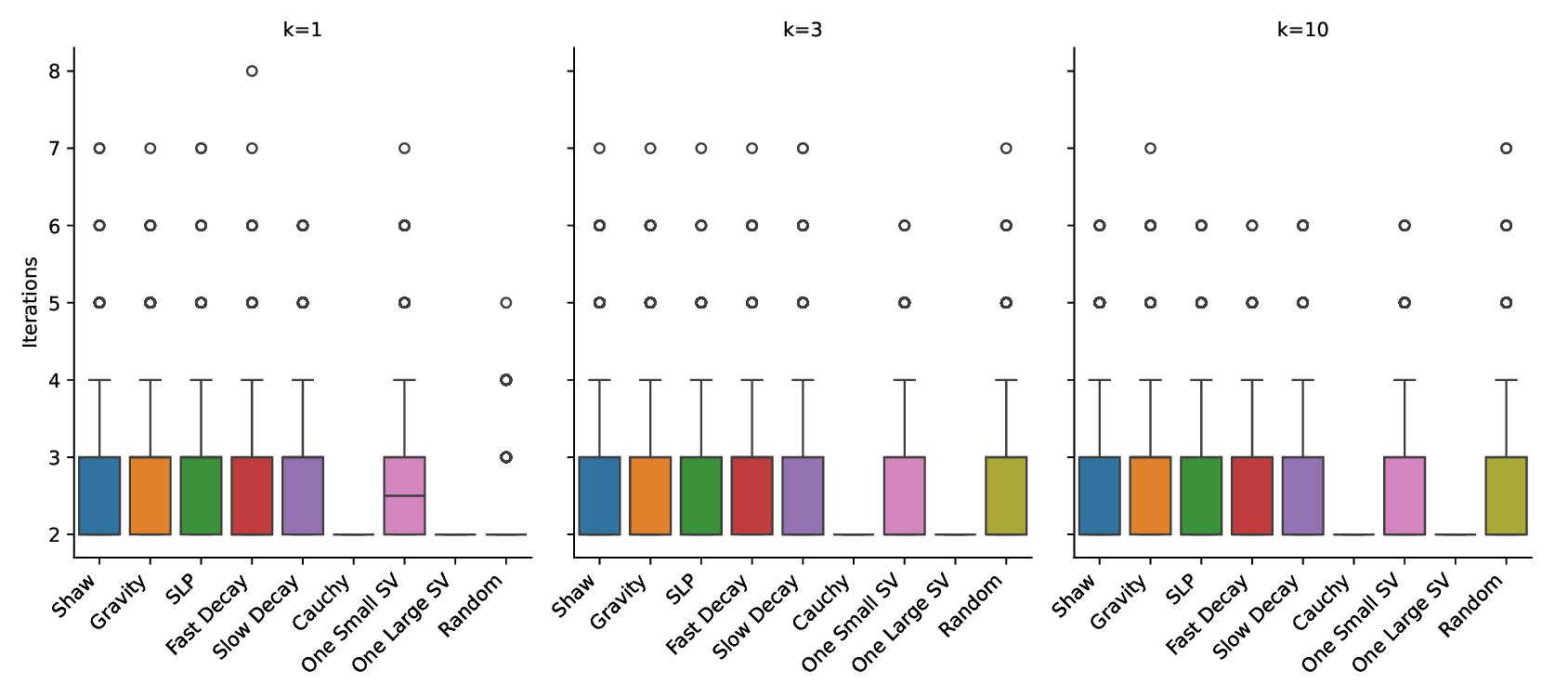}
\caption{Number of iterations of Alg.~\ref{algnormest}.}
\label{figalg22its}
\end{figure}

We ran 1000 tests for
Alg. \ref{algnormest}
for $\alpha = n/k$ under the same setup as for
Alg.~\ref{algnormest2}.  Figs. \ref{figalg22acc} and \ref{figalg22its} show the test results.  More refined data in Table~\ref{tabalg21-22comparison}   display   the mean output errors of both algorithms and the percent of the difference in average relative errors,
  $$\Delta {\rm err}
=:100 \times \frac{\text{err}_{\text{Alg}\ref{algnormest2}} - \text{err}_{\text{Alg}\ref{algnormest}}}{\text{err}_{\text{Alg}\ref{algnormest2}}}~{\rm where}~\text{err} := \text{mean} -1.$$ 
These data show  that both algorithms
yield similar output accuracy.

\begin{table}[h!]
\centering
\begin{tabular}{l||cc|c||cc|c||cc|c}
\toprule
 \multirow{2}{*}{Input Matrix}  
 & \multicolumn{3}{c||}{k=1} 
 & \multicolumn{3}{c||}{k=3} 
 & \multicolumn{3}{c}{k=10} \\
 & Alg.
 \ref{algnormest2}   & Alg.
 \ref{algnormest} & $\Delta$err\% 
 & Alg.
 \ref{algnormest2}   & Alg.
 \ref{algnormest} & $\Delta$err\% 
 & Alg.
 \ref{algnormest2}   & Alg.
 \ref{algnormest} & $\Delta$err\% \\
\midrule
Shaw & 1.1296 & 1.1407 & -8.54 & 1.0422 & 1.0438 & -3.93 & 1.0239 & 1.0276 & -15.80 \\
Gravity & 1.0536 & 1.0553 & -3.13 & 1.0300 & 1.0270 & 9.94 & 1.0248 & 1.0231 & 6.87 \\
SLP & 1.0013 & 1.0013 & -1.07 & 1.0009 & 1.0009 & -1.23 & 1.0003 & 1.0003 & -7.88 \\
Fast Decay & 1.1610 & 1.1622 & -0.71 & 1.1591 & 1.1531 & 3.77 & 1.1592 & 1.1647 & -3.48 \\
Slow Decay & 1.1540 & 1.1533 & 0.42 & 1.1618 & 1.1620 & -0.12 & 1.1596 & 1.1682 & -5.37 \\
Cauchy & 1.0000 & 1.0000 & - & 1.0000 & 1.0000 & - & 1.0000 & 1.0000 & - \\
One Small SV & 1.0222 & 1.0224 & -0.57 & 1.0212 & 1.0209 & 1.43 & 1.0206 & 1.0206 & -0.01 \\
One Large SV & 1.0000 & 1.0000 & - & 1.0000 & 1.0000 & - & 1.0000 & 1.0000 & - \\
Random & 1.0644 & 1.0645 & -0.14 & 1.0546 & 1.0541 & 0.82 & 1.0526 & 1.0526 & -0.09 \\
\bottomrule
\end{tabular}
\caption{Comparison of average output errors of Algs.~\ref{algnormest2} and \ref{algnormest}.}
\label{tabalg21-22comparison}
\end{table}

\subsection{Test Results for Algorithm \ref{algcrit}}\label{sslg3.2}

\begin{figure}[h!]
\centering
\begin{subfigure}{0.49\textwidth}
  \centering
  \includegraphics[width=.7\textwidth]{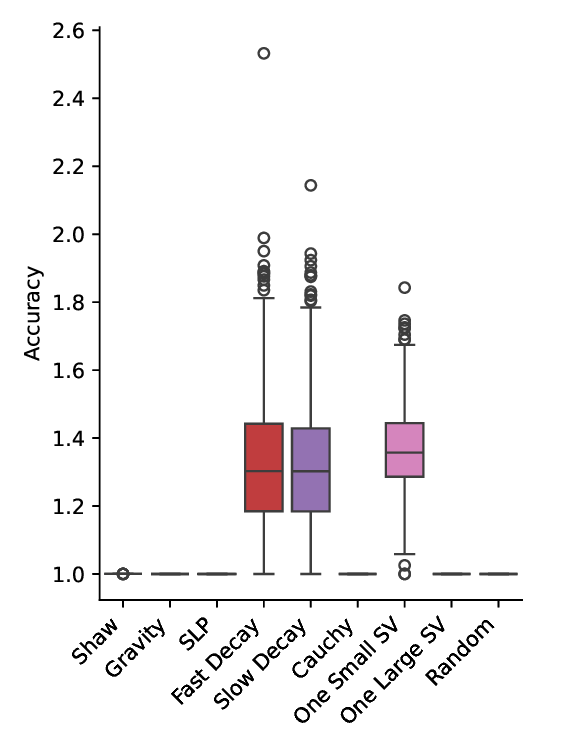}
  \caption{Output errors of Alg.~\ref{algcrit}}
  \label{figalg31acc}
\end{subfigure}
\hfill
\begin{subfigure}{0.49\textwidth}
  \centering
  \includegraphics[width=.7\textwidth]{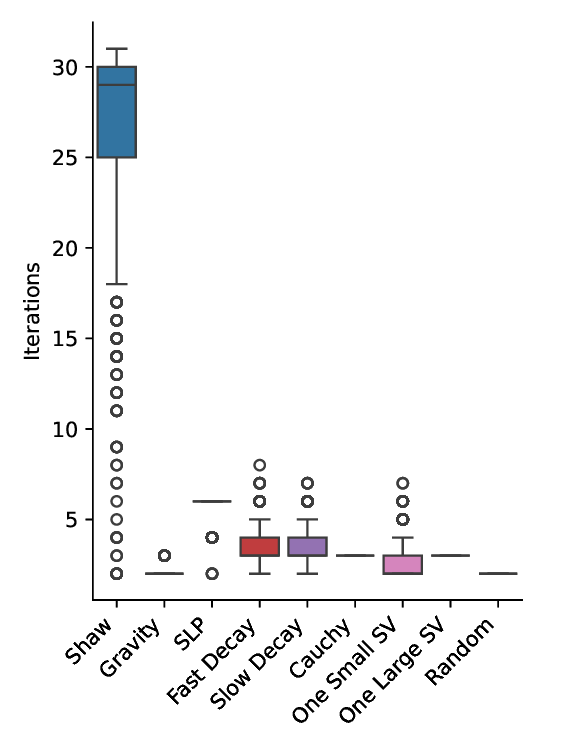}
  \caption{Number of C-A steps in  Alg.~\ref{algcrit}}
  \label{figalg31its}
\end{subfigure}
\caption{Output errors and the number of C-A steps of Alg.~\ref{algcrit}.}
\label{figalg31combined}
\end{figure}

We ran 1000 tests for Alg.~\ref{algcrit} for each type of input matrices. We measured output errors by the ratio $\frac{|M|}{\text{est}}$ where est denoted the estimated norm $|M|$
of (\ref{eqij||}) returned by the algorithm. Figs. \ref{figalg31combined}(a) and  \ref{figalg31combined}(b) show our error estimates and the numbers of argmaxima computed by the algorithm, including an argmaximum in the initialization step.  The  algorithm used
$\ll n^2$ comparisons;  the value $|m_{i,j}|$
for the output 
pair  $i$ and $j$ tended to be within a small factor from $|M|$.

\subsection{Test Results for Algorithm~\ref{algnormest1}}
\label{salg3.2}

\begin{figure}[h!]
\centering
  \includegraphics[width=.9\textwidth]{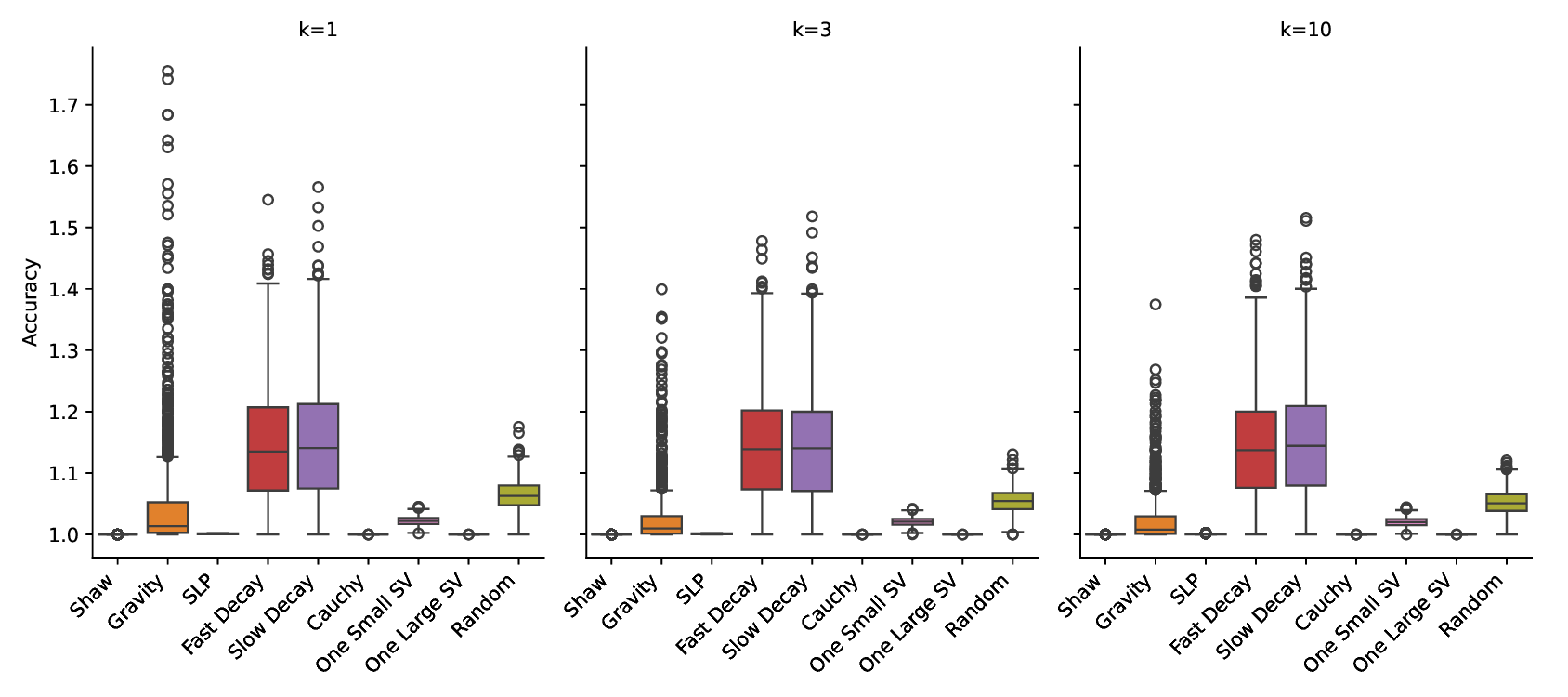}
  \caption{Accuracy (error estimates) of Alg.~\ref{algnormest1}.}
  \label{figalg32acc}
\end{figure}

\begin{figure}[h!]
\centering
  \includegraphics[width=.9\textwidth]{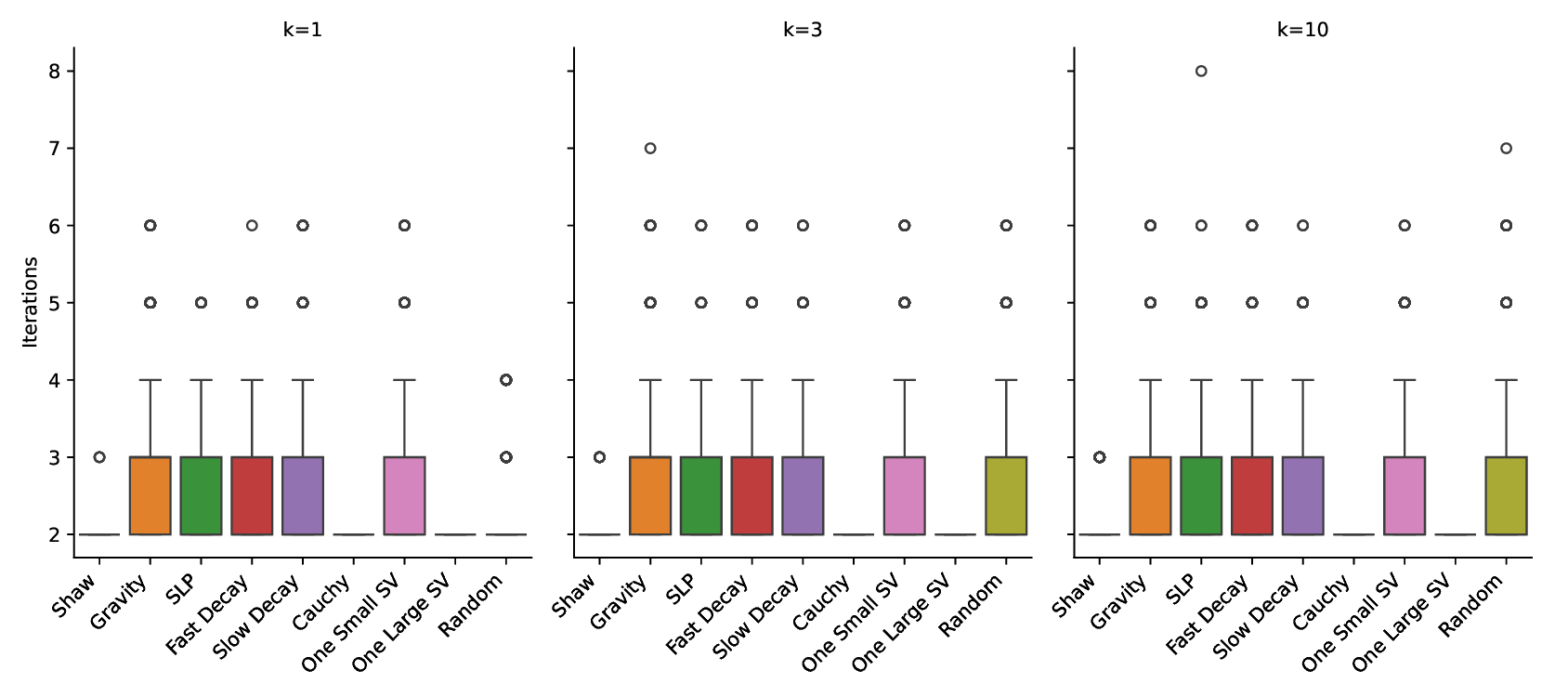}
  \caption{Number of Iterations of  Alg.~\ref{algnormest1}.}
  \label{figalg32its}
\end{figure} 

We ran 1000 tests for Alg.~\ref{algnormest1} under the same setup as for Alg.~\ref{algnormest2} 
with an additional input parameter $\text{tol}$, which we set to 1. 
Figs. \ref{figalg32acc} and \ref{figalg32its} display the test results.
They are similar to those for Alg.~\ref{algnormest2}.  More refined data in  Table~\ref{tabalg21-32comparison} display the output errors of both algorithms and the percent of the difference in average relative errors. These data show  a slightly better output accuracy
of Alg.~\ref{algnormest1}. 

\begin{table}[h!]
\centering
\begin{tabular}{l||cc|c||cc|c||cc|c}
\toprule
 \multirow{2}{*}{Input Matrix}  
 & \multicolumn{3}{c||}{k=1} 
 & \multicolumn{3}{c||}{k=3} 
 & \multicolumn{3}{c}{k=10} \\
  & Alg.
  \ref{algnormest2} & Alg.
  \ref{algnormest1} & $\Delta$err\% 
  & Alg.  
  \ref{algnormest2} & Alg.
  \ref{algnormest1}  & $\Delta$err\%  
  & Alg.
  \ref{algnormest2} & Alg. 
  \ref{algnormest1}& $\Delta$err\%  \\
\midrule
Shaw & 1.1296 & 1.0000 & 99.98 & 1.0422 & 1.0000 & 99.93 & 1.0239 & 1.0000 & 99.87 \\
Gravity & 1.0536 & 1.0508 & 5.22 & 1.0300 & 1.0282 & 5.85 & 1.0248 & 1.0247 & 0.49 \\
SLP & 1.0013 & 1.0012 & 6.88 & 1.0009 & 1.0009 & 1.09 & 1.0003 & 1.0004 & -18.64 \\
Fast Decay & 1.1610 & 1.1446 & 10.21 & 1.1591 & 1.1432 & 9.97 & 1.1592 & 1.1417 & 11.01 \\
Slow Decay & 1.1540 & 1.1478 & 4.00 & 1.1618 & 1.1434 & 11.37 & 1.1596 & 1.1484 & 7.01 \\
Cauchy & 1.0000 & 1.0000 & - & 1.0000 & 1.0000 & - & 1.0000 & 1.0000 & - \\
One Small SV & 1.0222 & 1.0218 & 2.02 & 1.0212 & 1.0207 & 2.53 & 1.0206 & 1.0201 & 2.41 \\
One Large SV & 1.0000 & 1.0000 & - & 1.0000 & 1.0000 & - & 1.0000 & 1.0000 & - \\
Random & 1.0644 & 1.0642 & 0.28 & 1.0546 & 1.0550 & -0.79 & 1.0526 & 1.0518 & 1.43 \\
\bottomrule
\end{tabular}
\caption{Comparison of the average output errors of Algs.~\ref{algnormest2} and \ref{algnormest1}.}
  \label{tabalg21-32comparison}
\end{table}

\subsection{Test Results for Algorithm~\ref{algcrit}-k}
We performed 1000 tests of Alg.~\ref{algcrit}-k for $k = 1, 3,10$ and  
in Table \ref{tabalg31-33comparison}  compare  average output errors of 
Algs.~\ref{algcrit} and \ref{algcrit}-k, both choosing the initial column index from the output of Alg. ~\ref{algnormest2}.
The table shows similar output accuracy of both algorithms. 

\begin{table}[h!]
\centering
\begin{tabular}{l||c||cc|cc|cc}
\toprule
 &
& \multicolumn{2}{c|} {$k=1$} 
& \multicolumn{2}{c|}{$k=3$} 
& \multicolumn{2}{c}{$k=10$}\\
Input Matrix & Alg.~\ref{algcrit} &  Alg.~\ref{algcrit}-k& $\Delta$err\% &  Alg.~\ref{algcrit}-k& $\Delta$err\% &   Alg.~\ref{algcrit}-k& $\Delta$err\% \\
\midrule
Shaw & 1.0001 & 1.0001 & 1.14 & 1.0001 & 2.49 & 1.0001 & 2.69 \\
Gravity & 1.0000 & 1.0000 & - & 1.0000 & - & 1.0000 & - \\
SLP & 1.0000 & 1.0000 & -5.49 & 1.0000 & -6.59 & 1.0000 & -4.81 \\
Fast Decay & 1.3228 & 1.2711 & 16.02 & 1.2652 & 17.83 & 1.2638 & 18.26 \\
Slow Decay & 1.3197 & 1.2644 & 17.31 & 1.2639 & 17.46 & 1.2663 & 16.70 \\
Cauchy & 1.0000 & 1.0000 & - & 1.0000 & - & 1.0000 & - \\
One Small SV & 1.3656 & 1.3805 & -4.08 & 1.3665 & -0.26 & 1.3695 & -1.07 \\
One Large SV & 1.0000 & 1.0000 & - & 1.0000 & - & 1.0000 & - \\
Random & 1.0000 & 1.0000 & - & 1.0000 & - & 1.0000 & - \\
\bottomrule
\end{tabular}
\caption{Comparison of the average output errors of Algs.~\ref{algcrit} and \ref{algcrit}-k.}
\label{tabalg31-33comparison}
\end{table}

\subsection{Test Results for Algorithm~\ref{algnrmest}}

We performed 1000 tests of Alg.~\ref{algnrmest} and display the results for $p =1$ and $r = 8, 16, 32, 64$. We computed rank-$r$ approximation $\{X, Y, Z\}$ of the input matrices   by using their random sketches $FM$ and $MH$ computed  either superfast with
3-Abridged SRHT multipliers
$F$ and $H$ 
or fast with Gaussian multipliers
$F$ and $H$
where  $F\in \mathbb R^{2r \times n}$ and $H\in \mathbb R^{n \times r}$
 (see Appendix~\ref{sspfcrd}). We  computed the 1-norm $||M||_1$ by using the LAPACK routine dlange accessed via scipy.linalg.lapack.dlange. In Table \ref{tabalgnrmestp1}  we display  the relative errors $\frac{||M||_1}{\text{est}}$ where est denoted the estimate of $||M||_1$ returned by the algorithm. 
The output accuracy of Alg. \ref{algnrmest} was noticeably inferior to that of Algs. \ref{algnormest2}  and \ref{algnormest1}
in our tests, but the relative errors  always stayed below 2.5 and  were similar for Gaussian and 3-Abridged  SRHT
 multipliers. 
\begin{table}[h!]
\centering
\begin{tabular}{l||c||cccc}
\toprule
 & Mutiplier & \multicolumn{4}{c} {$r$} \\
Input Matrix &  Type &  8 & 16 & 32 & 64  \\
\midrule
\multirow{2}{*}{Shaw}
 & 3-ASRHT & 1.1379 & 1.1594 & 1.1589 & 1.2168 \\
 & Gaussian & 1.0978 & 1.1579 & 1.2027 & 1.2358 \\\hline
\multirow{2}{*}{Gravity}
 & 3-ASRHT & 1.4349 & 1.4129 & 1.3930 & 1.3939 \\
 & Gaussian & 1.3019 & 1.3453 & 1.3765 & 1.4042 \\\hline
\multirow{2}{*}{SLP}
 & 3-ASRHT & 2.2551 & 1.9168 & 1.8366 & 1.7104 \\
 & Gaussian & 1.4883 & 1.5230 & 1.5560 & 1.5790 \\\hline
\multirow{2}{*}{Fast Decay}
 & 3-ASRHT & 0.9009 & 1.2754 & 1.3249 & 1.3447 \\
 & Gaussian  & 0.8895 & 1.2723 & 1.3278 & 1.3719 \\\hline
\multirow{2}{*}{Slow Decay}
 & 3-ASRHT & 0.9262 & 1.2745 & 1.3266 & 1.3668 \\
 & Gaussian & 0.9179 & 1.2743 & 1.3257 & 1.3639 \\\hline
\multirow{2}{*}{Cauchy}
 & 3-ASRHT & 1.8293 & 1.5966 & 1.3461 & 1.0539 \\
 & Gaussian & 0.6896 & 0.6567 & 0.6313 & 0.5967 \\\hline
\multirow{2}{*}{One Small SV}
 & 3-ASRHT & 0.4277 & 0.5088 & 0.6060 & 0.7114 \\
 & Gaussian & 0.4240 & 0.5085 & 0.6131 & 0.7256 \\\hline
\multirow{2}{*}{One Large SV}
 & 3-ASRHT & 1.4778 & 1.4595 & 1.4268 & 1.4259 \\
 & Gaussian & 1.5100 & 1.4657 & 1.4298 & 1.4254 \\\hline
\multirow{2}{*}{Random}
 & 3-ASRHT & 0.4683 & 0.5618 & 0.6682 & 0.7961 \\
 & Gaussian & 0.4561 & 0.5518 & 0.6668 & 0.8019 \\
\bottomrule
\end{tabular}
\caption{The average relative output errors of Algorithm \ref{algnrmest} for $p = 1$.}
\label{tabalgnrmestp1}
\end{table}

\section{Conclusions}\label{scnc}

We proposed, analyzed and tested a number of  superfast variations of 
 LAPACK's 1-norm estimator, which can be immediately extended to estimation of the $\infty$-norm. Our tests consistently show good output accuracy of
 our superfast  estimators.  Formal support of this empirical performance is an important challenge, possibly as hard as for C-A iterations.
Meanwhile further work  on pairwise comparison,  testing, and refinement
of our 1-norm estimators
are in order.

 One can further explore LRA-based $p$-norm  estimation for  $p\ge 1$ and $p=0$
(cf. Sec. \ref{sbckmnrm}) as well as the following approach to the same goal
 by means of randomized  sparsification with scaling:
(i) apply a known $p$-norm estimator
to a random $k\times \ell$ submatrix of an $m\times n$ input matrix $M$
for $k\ell\ll mn$,
(ii) repeat this computation $h$ times  for a fixed small $h$,
(iii)scale
the output estimates by a proper factor defined empirically,  repeat this computation $h$ times  for a fixed small $h$, and (iv)
output the average estimate.
 For the estimation of the Frobenius norm of $M$, one can modify this approach by keeping a  random subset of $k\ell$ entries of $M$ (rather than  its $k\times \ell$ submatrix) and setting to 0s the remaining entries.

Our main goal, however, was  to motivate more effort for
devising  superfast and empirically accurate  matrix algorithms, possibly by
means of extending our present techniques. E.g., one can
 explore  randomized sparsification with scaling for various matrix computations, but
 we conclude with  an immediate extension of 
superfast LRA to superfast  solution of the  
 Linear Least Squares Problem (LLSP):
  given $M\in \mathbb C^{m\times n}$
  and ${\bf b}\in \mathbb C^{m\times 1}$ for $m>n$,  approximate a vector
\begin{equation}\label{eqllsp}  
{\bf z}={\rm argmin}_{\bf y}||M{\bf y}-{\bf b}||_2.
\end{equation}
 
The authors of \cite{ALS24} point out 
the following  implication of the triangle inequality,
 $$||M{\bf y}-{\bf b}||_2\le ||\tilde M{\bf y}-{\bf b}||_2+||(\tilde M-M){\bf y}||_2,$$
 for an LRA 
 $\tilde M$
computed by means of C-A iterations
and then in \cite[Sec. 5]{ALS24} 
 show  extension to
the least squares approximation  of multivariate functions
by using a polynomial
basis. 

 Clearly, however, the latter study of \cite{ALS24}
can be extended 
to any  
LRA algorithm, e.g., the superfast ones of \cite{GLPa},
rather than just C-A algorithms.
Suppose that such a superfast algorithm has computed  two matrices $P$ and $S$ such that $\tilde M=PS\in \mathbb R^{m\times n}$, 
$P\in \mathbb R^{m\times r}$, 
$S\in \mathbb R^{r\times n}$, and $r\ll \min\{m,n\}$. Then we can first compute a solution
$${\bf x}={\rm argmin}_{\bf v}||P{\bf v}-{\bf b}||_2$$ of an LLSP for a matrix $P$ of   
 a small  size
$m\times r$, then immediately recover the solution
$\tilde {\bf z}$ of the LLSP
  $\tilde {\bf z}={\rm argmin}_{\tilde {\bf y}}||\tilde M\tilde {\bf y}-{\bf b}||_2$ 
from the
underdetermined 
linear system
${\bf x}=S\tilde {\bf z}$, and 
 notice that 
 $||(\tilde M-M){\bf y}||_2\le ||\tilde M-M||_2\cdot ||{\bf y}||_2\mapsto 0$
 as so does the norm  $||\tilde M
 -M||_2$.

\bigskip


{\bf \Large Appendix} 
\appendix


\section{Small families of matrices that are hard for fast LRA}\label{shrdin}

  Any algorithm that does not access all entries of an input matrix
  fails to estimate the norms as well as to compute a close LRA of  the following small families of matrices.
  
\begin{example}\label{exdlt} 
 Let  $\Delta_{i,j}$ denote an $m\times n$ matrix
 of rank 1  filled with 0s except for its $(i,j)$th entry filled with 1. The $mn$ such matrices $\{\Delta_{i,j}\}_{i,j=1}^{m,n}$ form a family of  $\delta$-{\em matrices}.
We also include the $m\times n$ null matrix $O_{m,n}$
filled with 0s  into this family.
Now fix any   algorithm that does not access the $(i,j)$th  
entry of its input matrices  for some pair of $i$ and $j$. Such an algorithm  approximates 
the matrices $\Delta_{i,j}$ and $O_{m,n}$,
with an undetected  error at least 1/2.
We arrive at the same conclusion by applying the same argument to the
set of $mn+1$ small-norm perturbations of 
the matrices of the above family. Likewise,
  we can verify that any
 randomized 
algorithm that does not access all entries of an input matrix
fails on those matrix families with 
error probability not close to 0.
\end{example}

\begin{remark}\label{redlt} 
The  matrices of the  above families  represent data singularities 
in contrast to matrices representing regular processes or, say, smooth surfaces as well as to various important special classes of structured matrices such as symmetric or diagonally dominant.  
\end{remark}


\section{LRA with dense and Abridged SRHT random  sketch matrices}
\label{sspfcrd}
                                                                                                                                                                                                                                                                                                                                                                                                                                                                                                                                                                                                                                                                                                                                                                                                                                                                                                                                                                                                                                                                                                                                                                                                                                                                                                                                                                                                                                                                                                                                                                                                                                                                                                                                                                                                                            
\subsection{LRA with dense   sketch matrices}
\label{sdrnd}

 
SVD of a matrix defines its optimal LRA
even under the spectral  norm, but randomized LRA algorithms in \cite{HMT11,M11,TYUC17,CW17,N20,TWa}, and the  references therein run faster
  than the known SVD algorithms
  \cite[Fig. 8.6.1]{GL13}
  and are still expected to output  accurate
  or even near-optimal LRAs.
  
 The generalized Nystr{\"o}m  algorithm of \cite{N20}  modifies the algorithms of
  \cite{CW09,TYUC17}  to
 improve numerical stability of rank-$r$ approximation of  
a matrix $M\in \mathbb C^{m\times n}$ for 
$r\le  n\le m$. It
first generates
two random   matrices $F\in \mathbb C^{k\times m}$ and $H\in \mathbb C^{n\times \ell}$, for $r\le k\le \ell\le n$, then  computes the {\em sketches} $FM$ and $MH$, and finally 
 computes rank-$r$ approximation
 $\{X,Y,Z\}$
  of $M$.
  
  In the best studied case  $F$ and $H$ are  Gaussian random
  matrices. Then for  
 $k=2r+1$ 
 and $\ell=2k$  the expected Frobenius error norm  of the error matrix $M-XYZ$ is within a factor of 2 from optimal -- it does not exceed the  optimal one for rank-$\frac{k}{2}$ approximation.                                                                                                                                                                                                                                                                                                                                                                                                                                                                                                                                                                                                                                                                                                                                                                                                                                                                                                                                                                                                                                                                   

The algorithm  involves $O((k+\ell)mn)$ 
 flops, $mn$ entries of $M$, and about 
 $km+n\ell+k\ell$ other scalars, which one can store by overwriting some  entries of $M$. 
 
With SRHT or SRFT (rather than Gaussian) multipliers $F$ and $H$, one can still obtain quite accurate  LRA (with a little higher error probability) by using
$O(mn\log(k\ell)))$ flops, $mn$ entries of $M$, and 
 $mk\log(k)+n\ell\log(\ell)$ other
 scalars, which can overwrite the $mn$ entries of $M$.
 
Such random sketching LRA algorithms are superfast except for the stage of the computation of the sketches $FM$  and  $MH$ and can be readily modified to
bound the relative error of the output LRA by any  $\epsilon>0$, with the estimated  sketch size increasing fast as  $\epsilon\mapsto 0$.  

\subsection{Abridged SRHT  matrices}\label{spreprmlt}

 
With  {\em sparse subspace embedding} \cite{C16,CDDRa,CFSa},
\cite[Sec. 3.3]{TYUC19},
\cite[Sec. 9]{MT20} one obtains significant acceleration although  not  superfast algorithms. 
According to \cite{L09},
 such acceleration tends to
 make the accuracy of the output LRAs somewhat less reliable. \cite{CFSa} partly overcomes
 this  problem for incoherent matrices.  Towards turning a matrix into incoherent one, we can  multiply it by SRHT or SRFT dense matrices
 \cite{CFSa}. This step is not superfast but
becomes superfast when we 
{\em abridge  SRHT multipliers},  sacrificing formal support for output accuracy. 
Namely, we
write: 
 

\begin{equation}\label{eqrfd}
H_{d,0}:=I_{n/2^d},~
H_{d,i+1}:=\begin{pmatrix}
H_{d,i} & H_{d,i} \\
H_{d,i} & -H_{d,i}
  \end{pmatrix}
  ~{\rm for}~i=0,1,\dots,d-1. 
\end{equation}
The $n\times n$ matrix $H_{d,d}$ is orthonormal 
up to scaling and
 has $2^d$ nonzero entries 
in every row and  column. 

Now define  $\ell\times n$ 
 matrix $R$ of uniform random  sampling of $\ell$ out of $n$ columns, the  $n\times n$ diagonal matrix 
  $D$  whose $n$ diagonal entries are independent random signs, i.e., random variables uniformly distributed on the pair
 $\{\pm 1\}$, and the $d$-{\em Abridged
 subsampled randomized Hadamard transform (SRHT)}
 matrix  $\sqrt{2^d/\ell~}RH_{d,d}D$, with
  $ 2^d$
nonzero entries
$\pm  \sqrt{2^d/\ell}$   per column. We can multiply $M\in \mathbb C^{m\times n}$  by that matrix by
involving $m2^d$ 
entries of $M$
and $(2^d-1)m$ flops.

 Likewise, by   sampling $k$ random rows of $H_{d,d}$ and multiplying 
 the resulting matrix by $\pm  \sqrt{2^d/\ell}~D$  obtain 
$d$-Abridged $k\times m$ SRHT matrix having
$2^d$ nonzero entries $\pm  \sqrt{2^d/k}$  per row. We can multiply it by $M$ by involving  $n2^d$ 
entries of $M$
and $(2^d-1)n$ flops.

E.g., for $d=3$
this means  $8n$  (resp. $8m$) entries of $M$ and $7n$
(resp. $7m$) flops.

A $t$-Abridged SRHT matrix   turns into
an SRHT matrix.

Table \ref{tabalgb1} displays the computational cost of randomized rank-$r$ approximation by means of the algorithm of \cite{N20}, when it uses Gaussian, SRHT and $d$-Abridged   
SRHT multipliers $F$ and $H$.

\begin{table}[h!]
\centering
\begin{tabular}{|l|c|c|c|cc|}
\toprule
Multipliers $F$ and $H$& Input entries involved & Other entries involved &  Flops involved  \\
\midrule
Gaussian & $mn$& $km+n\ell+k\ell$ & $O((k+\ell)mn)$ \\

SRHT& $mn$ & $O(mk\log(k)+n\ell\log(\ell))$ & $O(mn\log(k\ell))$\\
$d$-Abridged SRHT &$2^d(m+n)$ & $2^d(k+\ell)$ & $(2^d-1)(m+n)$ \\
\bottomrule
\end{tabular}
\caption{Computational cost of Algorithm  of \cite{N20} with 
Gaussian, SRHT, and $d$-Abridged SRHT multipliers $F$ and $H$.}\label{tabalgb1}
\end{table}
 

\section{Search for a scaling factor $\alpha$}\label{sscl}

The next algorithm    
recursively adjusts  scaling factor $\alpha$ of   Alg. \ref{algnormest}. Towards this goal we recursively invoke {\bf Alg. \ref{algnormest}a}, which
we obtain by  modifying  Alg. \ref{algnormest}. Namely, we write $\nu_0 := ||{\bf u}||_1$ at the Initialization stage 
and output the value $\gamma$ instead of $j$ and $||M{\bf e}_j||_1$. 

We begin our search for $\alpha$
by fixing the sparsification ratio $n/k$  as its initial value and then
recursively
 double or halve it
where Alg. \ref{algnormest}a outputs 1  (stopping prematurely) or TOL
(stopping  late), respectively. If neither occurs,  we stop and output the current value of $\alpha$. Furthermore,
if in its two successive applications  
Alg. \ref{algnormest}a  outputs 1
for $\alpha=\alpha_-$ followed by TOL for $\alpha=\alpha_+$  or vice versa, then we
  stop and output $\alpha:=\sqrt{\alpha_-\alpha_+}=\sqrt 2~\alpha_-$.
  
 Alg.~\ref{algalphasearch1} 
  implements this policy.
  
\begin{algorithm}\label{algalphasearch1}
{\rm [Search for $\alpha$ with doubling and halving.]}


\begin{description}


\item[{\sc Input:}] two  integers $k$ and $n$ such that $0<k\ll n$, a positive  TOL, and
a real $n\times n$ matrix $M$.


\item[{\sc Output:}]
A real $\alpha\ge 1$. 


\item[{\sc Initialization:}] \quad
 Set $\alpha \gets n/k$.  
 
 \item[{\sc Computations:}] \quad
Stop and output $\alpha$ unless Alg. \ref{algnormest}a outputs 1 or TOL.

Stop and output
$\alpha/\sqrt 2$ if in two successive applications it outputs 1 and then TOL.

Stop and output
$\alpha \sqrt 2$ if  in two successive applications it outputs  TOL and then 1. 

Double  $\alpha$ if
Alg. \ref{algnormest}a outputs 1.

Halve $\alpha$ if
it outputs TOL.

Repeat.

\end{description}
\end{algorithm} 

\begin{remark}\label{rebnrsrch} 
 We can further  narrow the range $(\alpha_-,\alpha_+)$ by applying binary search for 
  $\log_2(\alpha)$. At  every recursive step of the search we narrow the input range $[\alpha_-,\alpha_+]$
  to one of the two ranges $[\alpha_-,\beta]$ or $[\beta,\alpha_+]$ for $\beta:=\sqrt{\alpha_-\alpha_+}$. We stop
and output $\beta:=\sqrt{\alpha_-\alpha_+}\ge 1$
(for the current range $[\alpha_-,\alpha_+]$) wherever either Alg. \ref{algnormest}a  outputs a value strictly between 1 and TOL or  the ratio
 $\beta-1$ decreases below a fixed bound.
\end{remark}

\begin{remark}\label{remjrtvt} 
In more refined and more costly policy, we decide about doubling, halving or stopping according to majority vote in $h$ successive applications 
of Alg. \ref{algnormest}a for a fixed positive integer $h$. Namely, 
double $\alpha$
if more than 
 $h/2$ times the output COUNTER is 1;  halve
$\alpha$
if more than $h/2$ times the output COUNTER is TOL; otherwise  stop and output $\alpha$.
Hereafter refer to this modification of  
Alg.~\ref
{algalphasearch1} as {\bf Algorithm~ \ref{algalphasearch1}-h}, which
turns into 
Alg.~\ref{algalphasearch1}
for $h=1$. 
\end{remark} 

The overall cost of  
 execution of Alg.
\ref{algnormest}a,
invoked $hp$ times, dominates  the overall {\bf  computational complexity} of 
 Alg.~\ref{algalphasearch1}-h, which  therefore is much less than $n^2$ where $hps \ll n/k$. 
 
 The complexity  is  affected by
the value of $\alpha$ only through the impact of $\alpha$  on $p$, which we can share for all matrices  
 defined by similar data distributions.

In our tests  Algs. \ref{algalphasearch1} 
 and \ref{algalphasearch1}-h  slightly increased  output accuracy
 of Alg.  \ref{algnormest}
for some input classes by  slightly increasing the computational cost, but the adaptive processes of these algorithms have  
usually still ended with $\alpha\approx n/k$.



\end{document}